\input epsf

\magnification=1200


\hsize=125mm             
\vsize=195mm             
\parskip=0pt plus 1pt    
\clubpenalty=10000       
\widowpenalty=10000      
\frenchspacing           
\parindent=8mm           

\let\txtf=\textfont
\let\scrf=\scriptfont
\let\sscf=\scriptscriptfont
\font\frtnrm =cmr12 at 14pt
\font\tenrm  =cmr10
\font\ninerm =cmr9
\font\sevenrm=cmr7
\font\fiverm =cmr5

\txtf0=\tenrm
\scrf0=\sevenrm
\sscf0=\fiverm

\def\rm{\fam0 \tenrm}

\font\frtnmi =cmmi12 at 14pt
\font\tenmi  =cmmi10
\font\ninemi =cmmi9
\font\sevenmi=cmmi7
\font\fivemi =cmmi5

\txtf1=\tenmi
\scrf1=\sevenmi
\sscf1=\fivemi

 \def\oldstyle{\fam1 \tenmi}

\font\tensy  =cmsy10
\font\ninesy =cmsy9
\font\sevensy=cmsy7
\font\fivesy =cmsy5

\txtf2=\tensy
\scrf2=\sevensy
\sscf2=\fivesy

\font\tenit  =cmti10
\font\nineit =cmti9
\font\sevenit=cmti7
\font\fiveit =cmti7 at 5pt

\txtf\itfam=\tenit
\scrf\itfam=\sevenit
\sscf\itfam=\fiveit

\def\it{\fam\itfam\tenit}
\font\tenbf  =cmb10
\font\ninebf =cmb10 at  9pt
\font\sevenbf=cmb10 at  7pt
\font\fivebf =cmb10 at  5pt

\txtf\bffam=\tenbf
\scrf\bffam=\sevenbf
\sscf\bffam=\fivebf

\def\bf{\fam\bffam\tenbf}

\newfam\msbfam       
\font\tenmsb  =msbm10
\font\sevenmsb=msbm7
\font\fivemsb =msbm5

\txtf\msbfam=\tenmsb
\scrf\msbfam=\sevenmsb
\sscf\msbfam=\fivemsb

\def\msb{\fam\msbfam\tenmsb}
\def\Bbb#1{{\msb #1}}

\newfam\scfam
\font\tensc  =cmcsc10

\txtf\scfam=\tensc

\def\sc{\fam\scfam\tensc}


\def\frtnmath{%
\txtf0=\frtnrm        
\txtf1=\frtnmi         
}

\def\frtnpoint{%
\baselineskip=16.8pt plus.5pt minus.5pt%
\def\rm{\fam0 \frtnrm}%
\def\oldstyle{\fam1 \frtnmi}%
\everymath{\frtnmath}%
\everyhbox{\frtnrm}%
\frtnrm }


\def\ninemath{%
\txtf0=\ninerm        
\txtf1=\ninemi        
\txtf2=\ninesy        
\txtf\itfam=\nineit      
\txtf\bffam=\ninebf      
}

\def\ninepoint{%
\baselineskip=10.8pt plus.1pt minus.1pt%
\def\rm{\fam0 \ninerm}%
\def\oldstyle{\fam1 \ninemi}%
\def\it{\fam\itfam\nineit}%
\def\bf{\fam\bffam\ninebf}%
\everymath{\ninemath}%
\everyhbox{\ninerm}%
\ninerm }



\def\text#1{\hbox{\rm #1}}

\def\proj{\mathop{\rm proj}\nolimits}

\def\cite#1{{\uppercase{#1}}}
\def\ref#1{{\uppercase{#1}}}
\def\label#1{{\uppercase{#1}}}



\def\topmatter{\null\firstpagetrue\vskip\bigskipamount}
\def\endtopmatter{\vskip2\bigskipamount}

\def\title#1{%
\vbox{\raggedright\frtnpoint
\noindent #1\par}
\vskip 2\bigskipamount}        


\def\shorttitle#1{\rightheadtext={#1}}              

\newif\ifThanks
\global\Thanksfalse

\def\author#1{\begingroup\raggedright
\noindent{\sc #1\ifThanks$^*$\else\fi}\endgroup
\leftheadtext={#1}\vskip \bigskipamount}

\def\endabstract{\endgroup}

\long\def\abstract#1\endabstract{\par
\begingroup\ninepoint\narrower
\noindent{\sc Abstract.\enspace}#1%
\vskip\bigskipamount\endabstract}

\def\section#1#2{\bigbreak\bigskip\begingroup\raggedright
\noindent{\bf #1.\quad #2}\nobreak
\medskip\endgroup\noindent\ignorespaces}

\def\proclaim#1{\medbreak\noindent{\sc #1.\enspace}\begingroup
\it\ignorespaces}
\def\endproclaim{\endgroup\bigbreak}

\def\remark#1{\medbreak\noindent{\sc Remark \enspace}
\begingroup\ignorespaces}


\newdimen\EZ

\EZ=.5\parindent

\newbox\itembox

\newdimen\ITEM
\newdimen\ITEMORG
\newdimen\ITEMX
\newdimen\BUEXE

\def\iteml#1#2#3{\par\ITEM=#2\EZ\ITEMX=#1\EZ\BUEXE=\ITEM
\advance\BUEXE by-\ITEMX\hangindent\ITEM
\noindent\leavevmode\hskip\ITEM\llap{\hbox
to\BUEXE{#3\hfil$\,$}}%
\ignorespaces}


\newif\iffirstpage\newtoks\righthead
\newtoks\lefthead
\newtoks\rightheadtext
\newtoks\leftheadtext
\righthead={\ninepoint\rm\hfill{\the\rightheadtext}\hfill\llap{\folio}}
\lefthead={\ninepoint\rm\rlap{\folio}\hfill{\the\leftheadtext}\hfill}
\headline={\iffirstpage\hfill\else
\ifodd\pageno\the\righthead\else\the\lefthead\fi\fi}
\footline={\iffirstpage\hfill\global\firstpagefalse\else\hfill\fi}

\leftheadtext={}
\rightheadtext={}


\def\Refs{\bigbreak\bigskip\noindent{\bf References}\medskip
\begingroup\ninepoint\parindent=40pt}
\def\endRefs{\par\endgroup}
\def\endref{}

\def\ref{\par}
\def\key#1{\item{\hbox to 30pt{[#1]\hfill}}}


\def\cline#1{\leftline{\hfill#1\hfill}}

\def\bR{\Bbb R}
\def\bC{\Bbb C}

\def\bZ{\Bbb Z}

\hyphenation{mono-dro-my para-meter para-metrization para-metrize discri-mi-nant}
\input epsf

\topmatter
\title{Quadratic vanishing cycles, reduction curves and reduction of the
monodromy group of plane curve singularities}
\author{Norbert A'Campo}
\shorttitle{Vanishing cycles and reduction curves.}
\endtopmatter

\S 1. {\bf Introduction}
Let $f:\bC^2 \to \bC$ be a polynomial mapping with $f(0)=0$ having an
isolated singularity at $0 \in \bC^2$. Let $B \subset \bC^2$ 
be a Milnor ball for the
singularity. A disc $D \subset \bC$ with center $0 \in \bC$ is 
called {\it small relative to} $B$ if for all $s \in D$ the intersection of
$\partial{B}$ with $\{p\in \bC^2 \mid f(p)=s \}$ is a transversal intersection
of smooth differential manifolds. 

A {\it deformation of} $f$ is a polynomial mapping 
$f:\bC^2 \times \bC \to \bC,\,
(p,t) \mapsto f_t(p),$ with $f_0(p)=f(p)$.
A deformation of $f$ is {\it small relative to the Milnor ball}
$B$ {\it and the disk} $D$ if for all $s \in D$ and all 
$t \in \bC,\, |t| \leq 1,$ 
the intersection of
$\partial{B}$ with $\{p\in \bC^2 \mid f_t(p)=s \}$ 
is a transversal intersection
of smooth differential manifolds and moreover, the map $f_t$ is regular on
$f_t^{-1}(\partial{D})\cap B$.

Let $d \in D$ be the intersection point of $\partial{D}$ with the positive real
axis. The fibers $f_t^{-1}(d)\cap B, \,|t| \leq 1,$ are all 
diffeomorphic to the Milnor fiber $F:=f^{-1}(d)\cap B$, and the
diffeomorphism can be chosen to be unique up to isotopy
if one follows the family $(f_{st}^{-1}(d))_{s\in [0,1]}$.

A small deformation $f_t$ of $f$ is a {\it morsification of} $f$ if for all 
$t \in ]0,1],$ the map $f_t$ in $B\cap f_t^{-1}(D)$ has only 
singularities of morse type. {\it A real morsification} of a  singularity 
with a real defining equation $f$ is a morsification
$f_t$ with $f_t$ real, such that all critical points of the maps 
$f_t:B \to \bC$ are real and the level set 
$f_t^{-1}(0) \cap B,\,t\in ]0,1], $
contains all the saddle points of the 
restriction of $f_t$ to $D^2=B\cap \bR^2$. The level set 
$f_t^{-1}(0) \cap D^2$
is then a divide for the singularity of $f$, see [AC2,G-Z].

Let $f_t$ be a morsification, such that for some $\bar{t}\in ]0,1]$ 
the restriction to $B$ has $\mu(f)$ distinct critical values. 
A simply closed essential curve $c$ on the surface 
$f_{\bar{t}}^{-1}(d)\cap B$ is a {\it quadratic
vanishing cycle} if there exists a continuously differentiable path
$\gamma:[0,1] \to D$ such that $\gamma(0)=d$, 
$\gamma(u), u \in [0,1[$, 
are regular values of the restriction of $f_{\bar{t}}$ to $B$ and in
$f_{\bar{t}}^{-1}([0,1])$  the curve 
$c$ is null homotopic. A simply closed essential
curve $c$ on the Milnor fiber $F$ is a quadratic
vanishing cycle if there exists a morsification $f_t$ and 
$|\bar{t}|\leq 1, \bar{t}\not=0,$
such that $c$ is mapped to a quadratic vanishing cycle on 
$f_{\bar{t}}^{-1}(d)\cap B$ by
the natural diffeomorphism class. 

{\it The geometric monodromy group of the singularity of} $f$ is the image
$\Gamma_f$
of the geometric monodromy representation
in the mapping class group of the surface with boundary $(F,\partial{F})$ 
of the fundamental group of the complement (as germ) 
of the discriminant in the
unfolding of the singularity of $f$ [T2]. It follows from the results 
and constructions  of Egbert
Brieskorn [B] and a
theorem of Helmut
Hamm and 
L\^e D\~ung  Tr\'ang [H-L], that the geometric monodromy group 
is generated by the Dehn twists whose core curves are
the curves of any distinguished
system of vanishing cycles. 
Distinguished
systems of vanishing cycles for
isolated plane curve singularities have been
constructed by [AC2,G-Z].
The set of quadratic vanishing
cycles is an orbit of the action of the geometric monodromy group 
of the singularity. It follows that in principal the set of
quadratic vanishing cycles
is known.

A {\it reduction curve} for an isotopy class of  diffeomorphisms 
$T$ of a surface $F$
is an essential simply closed curve $a$ on $F$, such that 
for some integer $N > 0$ the curves 
$a$ and $T^N(a)$ are homotopic and moreover, if $N > 0$ is chosen minimal, 
the
curves $T^i(a), 0 \leq i \leq N-1,$ can be made pairwise disjoint by an 
isotopy of $T$. The set $\{T^i(a) \mid  0 \leq i \leq N-1\}$ is 
then called a reduction
system for $T$.

The set of reduction curves for the geometric monodromy of a
singularity with one local branch is given by [AC1]. The purpose of the present
paper is to describe the relative position of a distinguished 
base of vanishing
cycles and the set of reduction curves for the monodromy 
diffeomorphism of isolated plane
curve singularities. 

Roughly speaking, we wish to give a picture
of the Milnor fiber of a plane curve singularity, 
that shows both a distinguished system of quadratic vanishing cycles 
and the reduction curves of the geometric monodromy. The 
method uses the model of the Milnor fiber with 
monodromy coming from a divide of the singularity [AC2-5,G-Z]. It 
turns out that the
divides constructed by Sabir Gusein-Zade are most suited for this purpose.

In Section $4$ we will give an application of our constructions to 
the study of the geometric monodromy group of an irreducible plane curve
singularity with several essential Puiseux pairs.
It is proved that the geometric monodromy group of an irreducible
plane curve singularity in a natural way contains a product of many 
copies of  the geometric monodromy groups 
of its companion singularities.

\bigskip
\S 2. {\bf Building blocks and cabling} 

It was observed by Ren\'e Thom that the Chebyshev  polynomials $T:\bC \to \bC$ 
up to affine equivalence of functions are precisely the  polynomial mappings
from $\bC$ to $\bC$ with two or less critical values and with only quadratic
singularities [T1]. The standard Chebyshev 
polynomial $T(p,z)$ of degree $p$ has
critical values $+1,-1$, the symmetry $T(p,z)=(-1)^p T(p,-z)$, and the
coefficient of $z^p$ is $2^{p-1}$. For $p=1$, the map $T(p,z)$ has 
no critical values
and for $p=2$, the map $T(p,z)$ has only the critical value $-1$. The 
Chebyshev  polynomial $T(p,z)$ satisfies the identity
$T(p,\cos(x))=\cos(px)$, and its restriction to $[-1,1]$ is defined by
$T(p,t)=\cos(p\arccos(t))$.

Sabir Gusein-Zade has constructed real morsifications for
real plane curve singularities with Chebyshev
polynomials. The building blocks for his construction 
are the real morsification for the map $f(x,y)=2^{p-1}x^p-2^{q-1}y^q$ 
given by the family 
$$
f_s(x,y)=s^{pq}(
T(p,x/s^q)- T(q,y/s^p)), s \in [0,1].
$$
For each $s \in ]0,1]$ the
function $f_s:\bC^2 \to \bC$ has $\mu_f=(p-1)(q-1)$ quadratic 
singularities all at points with real coordinates, its critical values  
are contained in $\{-2s^{pq},0,2s^{pq}\}$
and $\lim_{s \to 0,s>0}f_s=f$.

If the exponents $p$ and $q$ are relatively prime to each other, the 
level set $\{f(x,y)=0\}$ can be
parametrized by the monomial map 
$$
t\in \bC \mapsto (t^q/2^{p-1},t^p/2^{q-1}) \in \bC^2.
$$
It is a miracle that in this case 
the Chebyshev polynomials can be used to parametrize the 
level sets $\{f_s(x,y)=0\}$ as well. The map 
$$
t \in \bC \mapsto (s^q T(q,t/s),s^p T(p,t/s)) \in \bC^2
$$ 
parametrizes the level set $\{f_s(x,y)=0\}$.
We recall that a divide is the image of a generic relative immersion
of a compact $1$-manifold in the unit disk $D^2$ in $\bR^2$.
For each $s\in ]0,1]$ the intersection $P_{p,q;s}:=\{f_s(x,y)=0\}\cap D^2$
with the unit disk $D^2$
is a divide
for the singularity of $\{2^{p-1}x^p-2^{q-1}y^q=0\}$ at $0\in \bC^2$. 

\midinsert
\cline{\epsffile {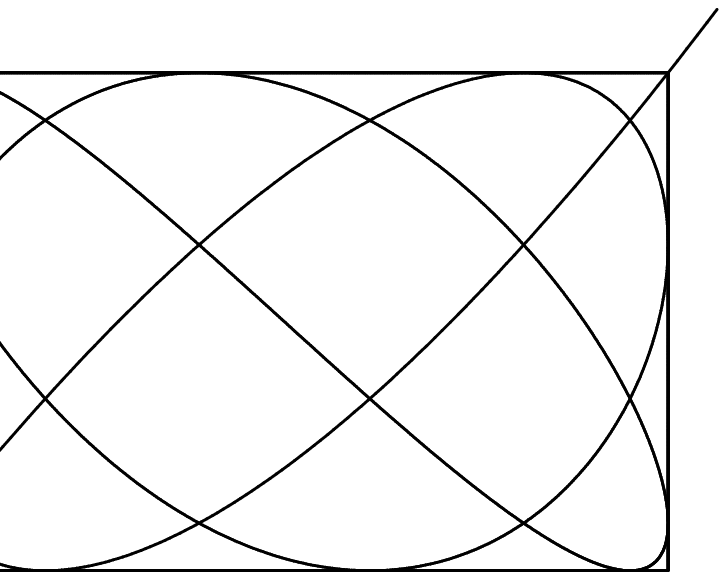}}
\centerline{Fig. $1.$ Divide in box $[-1,1]\times [-1,1]$ for  $2^6x^7-2^4y^5=0$.}
\endinsert

The curve  
$P_{p,q}:=\{f_1(x,y)=0\}$ can be drawn in a rectangular box as in Fig. $1$. 
As a first type of building block 
we will need 
the box $B:=[-1,1]\times [-1,1]$ with the curve $P_{p,q}$. 
If $(p,q)=1$ 
holds, the curve $P_{p,q}$ is the image of
$$
T_{p,q}:[-1,1] \to B, T_{p,q}(t):=(T(p,t),T(q,t)),
$$
which leaves the box through the corners. In general the immersed 
curve has several components, which are immersions of the interval or of the
circle. At most two components are immersions of the interval, 
which leave the box through the corners. 

Let $P$ be any divide having one branch given by an immersion 
$\gamma: [-1,1] \to D^2$. We assume, that the speed vector 
$\dot{\gamma}(t)$ and the position vector $\gamma(t)$ 
are proportional at $t=\pm 1$, i.e. the
divide $P$ meets $\partial{D^2}$ at right angles.
Let $N\gamma:[-1,1]\times [-1,1] \to D$ 
be the corresponding immersion of a
rectangular box, i.e. the restriction of 
$N\gamma$ to $[-1,1]\times \{0\}$
is the immersion $\gamma$ and the image of $N\gamma$ is in a small tubular
neighborhood of $P$. For instance, 
for a small value of the parameter $\eta \in \bR_{>0}$
the following
expression defines
an immersion $N\gamma:B \to \bR^2$ of the  rectangular box 
$B:=[-1,1]\times [-1,1]$: 
$$
N\gamma(s,t):=\gamma(t)+s\eta {J(\dot{\gamma}(t)) \over ||\dot{\gamma}(t)||},
$$
where $J$ is the rotation of $\bR^2$ over $\pi/2$.
The four corners $N\gamma(\pm 1,\pm 1)$ are on the circle of
radius 
$\sigma:=\sqrt{1+\eta^2}.$ 
We finally define 
$$
N\gamma(s,t):=N_{\eta}\gamma(s,t):={1 \over \sigma}(\gamma(t)+
s\eta {J(\dot{\gamma}(t)) \over ||\dot{\gamma}(t)||}),
$$
that will be an immersion $N\gamma:B \to D^2$ mapping the corners of 
the box $B$ 
into $\partial{D^2}$. 

We will denote by $P_{p,q}*P$ the divide in $D^2$, 
which is the image by $N\gamma:B \to D^2$ of 
$P_{p,q} \subset B$, see Fig. $2$.
The number of double
points $\delta(P_{p,q}*P)$ of $P_{p,q}*P$ is computed inductively 
from the number of double points $\delta(P)$ of $P$ by:
$$
\delta(P_{p,q}*P)=(p-1)(q-1)/2+\delta(P)p^2.
$$

Let $R_{\eta}\gamma$ be the union of the image of 
$N_{\eta}\gamma$ with the two
chordal caps at the endpoints of $\gamma$. The connected components
of $D^2 \setminus R_{\eta}\gamma$ correspond via inclusion to
the connected components of $D^2 \setminus P$. We declare a connected 
component of $D^2 \setminus P_{p,q}*P$ to be signed by $+$, if the component 
contains a component of $D^2 \setminus R_{\eta}\gamma$, that corresponds
to a $+$ component of $D^2 \setminus P$. In this case we 
call the connected component
of $D^2 \setminus R_{\eta}\gamma$ a $P_+$-{\it component}. Observe 
that there exists
a chess board sign distribution for the components of 
$D^2 \setminus R_{\eta}\gamma$ that makes $P_+$-components indeed to $+$
components.

The field $\Phi_{p,q}$ of cones on the box $B \subset \bR^2$ is 
the subset in the tangent space of 
$TB$ given by:
$$
\Phi_{p,q}:=\{(x,u) \in TB \mid |<u,e_1>_{\bR^2}| \geq \cos(\alpha(x))||u||\}
$$
where $e_1=(1,0)\in \bR^2$ and where $\alpha:B \to \bR$ is a function, such
that for every $(x,u) \in TB$ with $x \in P_{p,q}$ and $u \in T_xP_{p,q}$
we have the equality 
$$
|<u,e_1>_{\bR^2}| = \cos(\alpha(x))||u||.
$$
Moreover,
$\alpha$ has the boundary values 
$\alpha(\pm 1,t)=0$ and $\alpha(s,\pm 1)=\pi/2$.
We interpolate
the function $\alpha$ on $B$ by upper and lower convexity, i.e
such that ${\partial^2\over{\partial{t}}}\alpha< 0$ and 
${\partial^2\over{\partial{s}}}\alpha> 0$.
The definition of $\alpha(x)$ seems to be 
conflicting at the double points of
the curve $P_{p,q}$; 
at a double point $x=(x_1,x_2)$ of the curve 
$P_{p,q}$ the two tangents lines to $P_{p,q}$
have opposite slopes $\tan(\alpha(x))$ and $-\tan(\alpha(x))$, since
the curve $P_{p,q}$ is defined by  the  equation 
$$
T(q,x_1)-T(p,x_2)=0
$$ 
that separates the variables. 
For example, a nice such function
$\alpha$ is given by:
$$
\alpha(x_1,x_2):=\arctan({q\sqrt{1-x_2^2} \over p\sqrt{1-x_1^2}}).
$$

The interest of the field $\Phi_{p,q}$ comes from the following 
lemma, that is
immediate from the definitions.

\proclaim{Lemma} Let the image of $\gamma:[-1,1] \to D$ be a divide $P$, that
meets $\partial{D^2}$ at right angles.
For $\eta > 0$ small enough, the intersection of $S^3 \subset T\bR^2$ 
with the image in $T\bR^2$ 
of the
field of sectors $\Phi_{p,q}$ under the differential 
of $N_{\eta}\gamma$ is a
tubular neighborhood of the knot $L(P)$.
The composition of $P_{p,q}:[-1,1]  \to B:=[-1,1] \times [-1,1]$ and of
$N\gamma:B \to D^2$ 
is again a divide, whose knot is a torus cable knot of the knot $L(P)$.  
\endproclaim

The image of the field of sectors $\Phi_{p,q}$ of $B$ under 
the differential of 
$N_{\eta}\gamma$ will be denoted by 
$\Phi_{\eta,p,q}\gamma$ and for small $\eta$ contains 
those vectors, that 
have feet near $P$ and form a small 
angle with the tangent vectors of the divide $P$.

\midinsert
\cline{\epsffile {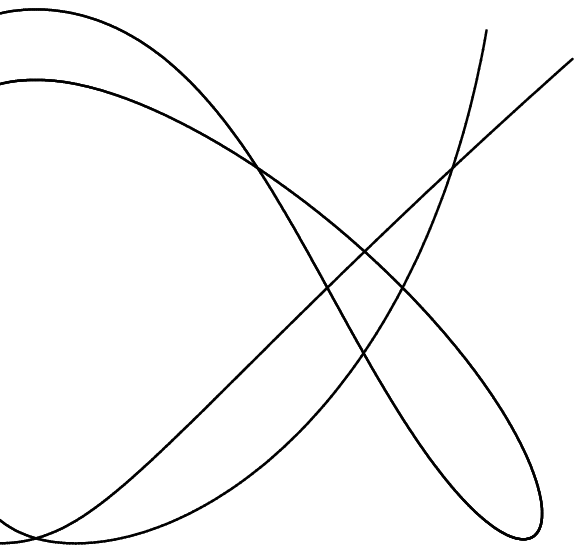}}
\centerline{Fig. $2.$ The divide $P_{2,9}*P_{2,3}$.}
\endinsert

A sequence of pairs of integers
$(a_i,b_i)_{1\leq i \leq k}$ is a {\it sequence of essential Puiseux pairs of
an irreducible plane curve singularity} if the inequalities
$2\leq a_i < b_i$ and
$b_i/a_1a_2 \cdots a_i < b_{i+1}/a_1a_2 \cdots a_ia_{i+1}$
are verified and if moreover, the integers $b_i$ and $a_1a_2 \cdots a_i$ 
are relatively prime. A sequence of essential Puiseux pairs defines a family
of topologically equivalent singularities. A specific member $f_k(x,y)$
of this family is
obtained from the Puiseux expansion with fractional and strictly increasing
exponents
$$
y=x^{b_1/a_1}+x^{b_2/a_1a_2}+ \dots +x^{b_k/a_1a_2 \cdots a_k}
$$
by the rule, which takes in account the ramification of $x^{1/a_1a_2\cdots a_k}$,
$$
f_k(x,y)=\prod_{\theta}(y-\theta^{a_2 \cdots a_kb_1}-\theta^{a_3 \cdots a_kb_2}- 
\dots -\theta^{b_k}),
$$
\noindent
where $\theta$ runs 
over the 
$a_1a_2 \cdots a_k$ roots 
of $z^{a_1a_2 \dots a_k}-x=0$ in the algebraic
closure of the field $\bC((x))$. The 
coefficients of the polynomial $f_k(x,y)$ are integers. 

For example,
the Puiseux 
expansion $y=x^{3/2}+x^{7/4}$ 
leads to the polynomial $f_{(2,3),(4,7)}=(y^2-x^3)^2-4x^5y-x^7$ and
the Puiseux expansion $y=x^{3/2}+x^{11/6}$ to the polynomial
$(y^2-x^3)^3-6x^7y^2-2x^{10}-x^{11}$.

Let $\{f_{a,b}(x,y)=0\}$ be a
singularity having one branch and with essential Puiseux pairs
$(a_i,b_i)_{1\leq i \leq n}$.
The theorem of S. Gusein-Zade [G-Z] very efficiently
describes a divide
for the singularity $\{f_{a,b}(x,y)=0\}$ in a closed form, namely
the iteratively composed divide
$$
P_{a_n,b'_n}* \dots *P_{a_2,b'_2}*P_{a_1,b_1},
$$
where the numbers $b'_2, \dots ,b'_n$ can be computed recursively, as
we will show here below.
We denote by $S_k, 1\leq k \leq n$ the divide
$$
P_{a_k,b'_k}*P_{a_{k-1},b'_{k-1}} \dots *P_{a_2,b'_2}*P_{a_1,b_1}
$$ and
let $f_k(x,y)$ be a specific  equation
for a singularity with essential Puiseux pairs $(a_i,b_i)_{1\leq i \leq k}$.

Remember, that the product $a_1 a_2 \dots a_k$ is the
multiplicity at $0$ of the curve $\{f_k(x,y)=0\}$ and that 
the linking number $\lambda_k$ of 
$L(S_k)$ and $L(S_{k-1})$ in $S^3$ can be computed recursively by:
$$
\lambda_1=b_1,\, \lambda_{k+1}=b_{k+1}-b_ka_{k+1}+\lambda_ka_ka_{k+1}.
$$
The linking number $\lambda_k$ is equal to the intersection multiplicity 
$$
\dim \bC[[x,y]]/(f_k(x,y),f_{k-1}(x,y))
$$ 
at $0$
of the curves $\{f_k(x,y)=0\}$ and $\{f_{k-1}(x,y)=0\}$.
We have also for the linking number $\lambda_k$ an interpretation in terms 
of divides (see the next section for the first equality)
$$
\lambda_k = \#(S_k \cap S_{k-1})=2a_k \delta(S_{k-1}) + b_k'.
$$
Remembering, that we already have 
computed recursively the numbers $\delta(S_k)$ and $\#(S_k \cap S_{k-1})$,
we conclude that $b_k'$ too can be computed recursively. 

\midinsert
\cline{\epsffile {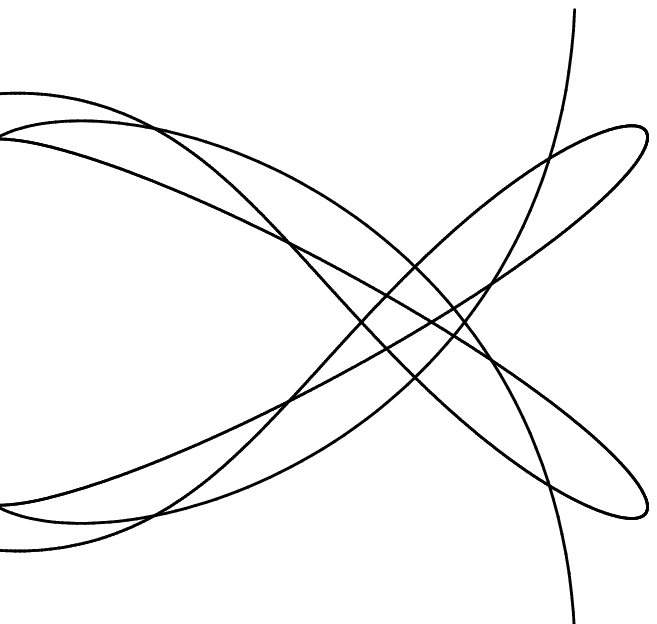}}
\centerline{Fig. $3.$ The divide $P_{3,14}*P_{2,3}$ for
$(y^2-x^3)^3-6x^7y^2-2x^{10}-x^{11}$.}
\endinsert

For  example, for
the Puiseux expansion $y=x^{3/2}+x^{7/4}$ we 
have: $\lambda_1=3,\, \lambda_2=7-3.2+3.2.2=13,\, b'_2=13-2.2.1=9$.
Hence,  the divide for the irreducible singularity with Puiseux expansion
$y=x^{3/2}+x^{7/4}$ is the divide $P_{2,9}*P_{2,3}$, see Fig. $2$. 
For the Puiseux expansion
$y=x^{3/2}+x^{11/6}$ we found: $\lambda_1=3,\, \lambda_2=11-3.3+3.2.3=20,\,
b'_2=20-2.3.1=14$. Hence, the divide for 
its singularity 
$\{(y^2-x^3)^3-6x^7y^2-2x^{10}-x^{11}=0\}$ is 
$P_{3,14}*P_{2,3}$, 
see Fig. $3$.

An iterated $*$-composition of divides
has to be evaluated from
the right to the left.

Using [AC4-5], we can read off from this divide the Milnor fibration of the
singularity $\{f_{a,b}(x,y)=0\}$. In particular we can describe the Milnor
fiber with a distinguished base of quadratic vanishing 
cycles, see next section. Using the above
iterated cabling construction, we will in Section $4$ also read off 
from the divide the reduction of
the geometric monodromy of an irreducible plane curve singularity, as 
described in [AC1]. For instance, 
intersection numbers in the sense of Nielsen of quadratic vanishing cycles and
reduction cycles can be computed.

In general, for an isolated singularity of a real polynomial 
$f(x,y)$ having
several local branches, the divide $\{f_1(x,y)=0\} \cap D$ of a 
real morsification $f_t(x,y)$ may have immersed circles as 
componants. The above cabling construction 
$P_{p,q}*P$ does not work if the divide $P$ consists
of an immersed circle. Of course, if one is willing to 
change the equation
of the singularity to an equation, which defines a 
topologically equivalent
singularity and which has only real local branches, one will only
have to deal with divides consisting of 
immersed intervals. If we do not want
to change the real equation, we will need a 
second type of building blocks for
a cabling construction, see Fig. $4$.

\midinsert
\cline{\epsffile {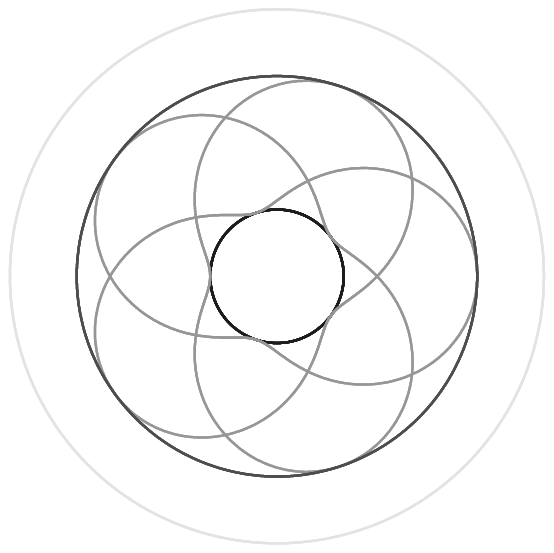}}
\centerline{Fig. $4.$ The building block  $L_{3,5}.$ of Lissajous type}
\endinsert

These building blocks are the divides 
$L_{p,q}$
in the annular region 
$A:=\{(x,y) \in D^2 \mid 1/4 \leq \sqrt{x^2+y^2} \leq 3/4 \}$.
If for the integers $(p,q)=1$ holds, the divide $L_{p,q}$ is the Lissajous
curve 
$$
s \in [0,1] \mapsto (1/2+1/4 \sin(2\pi q s))( \sin(2\pi p s),\cos(2 \pi p s))
$$
in $A \subset D^2$. The curve $L_{p,q}$ has $p$-fold rotational
symmetry. If $(p,q)=r > 1$ the divide $L_{p,q}$ will be 
defined as the union of 
$r$ rotated copies of 
$L_{p/r,q/r}$ with rotations of angles $2\pi k/p, k=0 \dots r-1,$ of $D^2$.
Again, the system of curves $L_{p,q}$ has a $q$-fold rotational symmetry.

The star-product $L_{p,q}*P$ can be defined as above 
if the divide $P$ consists
of one immersed circle. We leave many details to the reader. 
The two types of building blocks $P_{p,q}$ and $L_{p,q}$ 
together with the star-products $P_{p,q}*P$ and $L_{p,q}*P$ will 
allow one to describe the iterated cablings of real plane curve
singularities in general. See Fig. $5$.

\midinsert
\cline{\epsffile {lissa_2.4_3.5.eps}}
\centerline{Fig. $5.$ The cabling of  $L_{3,5}*P_{2,4}.$}
\endinsert

\bigskip
\S 3.  {\bf Visualization of the vanishing cycles for a divide}\par

Let $P$ be a connected divide and let 
$\pi_P:S^3 \setminus L(P) \to S^1$ be the 
fibration of the complement of the link $L(P)$ over $S^1$ as in [AC4-5].
The fibration map is given with the help of an auxiliary morse function
$f_P:D^2 \to \bR$.
The fiber $F_P$ above $1 \in S^1$ projects to the positive components of
the complement of $P$ in $D^2$. One has that the closure of 
$$
\{(x,u)\in TD^2 \mid f_P(x)>0, (df_P)_x(u)=0, ||x||^2+||u||^2=1\}
$$
in $S^3 \setminus L(P)$ is the fiber surface $F_P$.

To each critical point of $f_P:D^2 \to \bR$ corresponds a vanishing cycle on
the surface $F_P$. In the case, where the divide $P$ is a divide of a 
singularity, the surface $F_P$ is a model for the Milnor fiber and the
system of
vanishing cycles on $F_P$ is a model for  a distinguished system of quadratic
vanishing cycles of the singularity.

Let $M$ be a maximum of $f_P$. The vanishing cycle $\delta_M$ is the 
non-oriented simply closed curve on $F_P$, see Fig. $6$,
$$
\delta_M:=\{u\in T_MD^2 \mid ||M||^2+||u||^2=1\}. 
$$

\midinsert
\cline{\epsffile {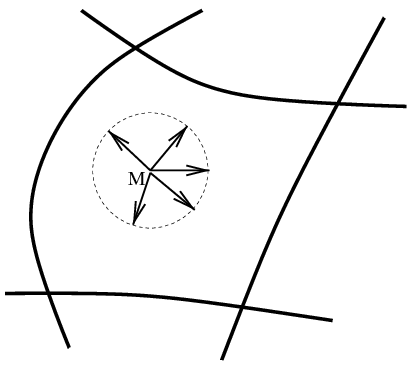}}
\centerline{Fig. $6.$ The vanishing cycle $\delta_M$ for a maximum.}
\endinsert

Let $c$ be a crossing point of $P$. The point $c$ is a saddle type 
singularity of $F_P$. The vanishing cycle is the non-oriented simply closed
curve $\delta_c$ on $F_P$ that results from the following construction.
Let $g_c \subset P_+ \cup \{c\}$ be the singular 
gradient line through $c$, for which
the endpoints are a maximum of $f_P$ or a point in $\partial{D^2}$. We splice
$g_c$ and get a double tear $t_c \subset P_+ \cup \{c\}$ as in Fig. $7$. The
tear $t_c$ is a closed curve, that has at $c$ a non-degenerate
tangency with $g_c$ from both sides. Moreover, $t_c$ 
is perpendicular to $g_c$ at the
endpoints of $g_c$, if the endpoint is a maximum of $f_P$ and else $t_c$
has a tangency with $\partial{P}$. The vanishing cycle $\delta_c$
is the closure in $F_P$ of the set 
$$
\{(x,u)\in F_P \mid  x\in t_c,\,\,(df_P)_x\not=0,\,u\,
\hbox{{\it points to the inside of the tear}}\,\,t_c\}.
$$

\midinsert
\cline{\epsffile {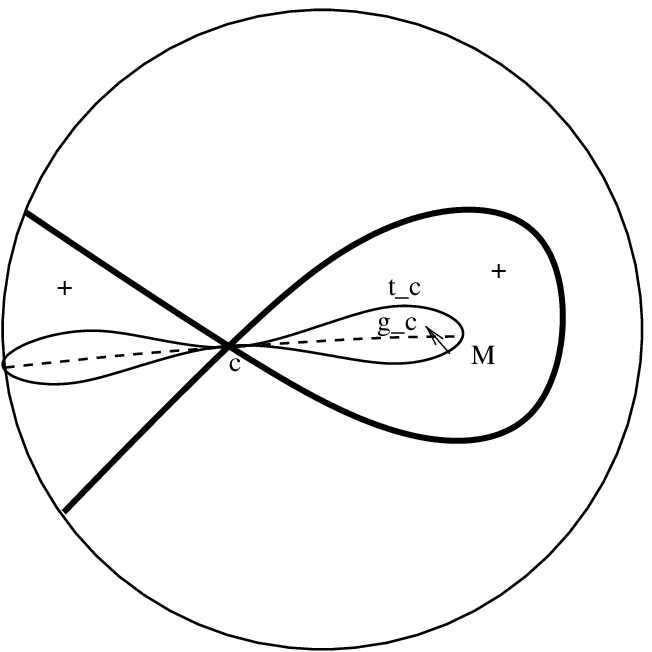}}
\centerline{Fig. $7.$ Tear for vanishing cycle $\delta_c$ for a saddle point.}
\endinsert 

Let $m$ be a minimum of $f_P$. The following is a description of the
vanishing cycle $\delta_m$ on $F_P$. The projection of $\delta_m$ in 
$D^2$ is a simply closed curve $t_m$ in $P_+ \cup \{c_1,c_2, \dots ,c_k\}$,
where $\{c_1,c_2, \dots ,c_k\}$ is the list of the double points of $P$
that lie in the closure of the region of $m$, see Fig. $8$. The curve 
$t_m$ and the singular gradient line $g_{c_i},\, 1\leq i\leq k$ coincide in
a neighborhood of $c_i$. Moreover, if the endpoint of $g_{c_i}$ is a 
maximum $M$,
the curve $t_m$ leaves transversally the tear $t_{c_i}$ at $M$ and enters
transversally the next tear $t_{c_{i\pm 1}}$. If
the endpoint of $g_{c_i}$ is on the boundary of $D$, the curve leaves the tear
$t_{c_i}$, becomes tangent to the boundary of $D^2$ and 
enters the next tear, see Fig. $8$. The vanishing cycle $\delta_m$ is the 
non-oriented simply closed curve on $F_P$, that is the closure of
$$\{(x,u)\in F_P \mid x \in t_m,\,u\,\hbox{{\it points inwards to the
disk bounded by }}\,\,t_m\}.
$$

\midinsert
\cline{\epsffile {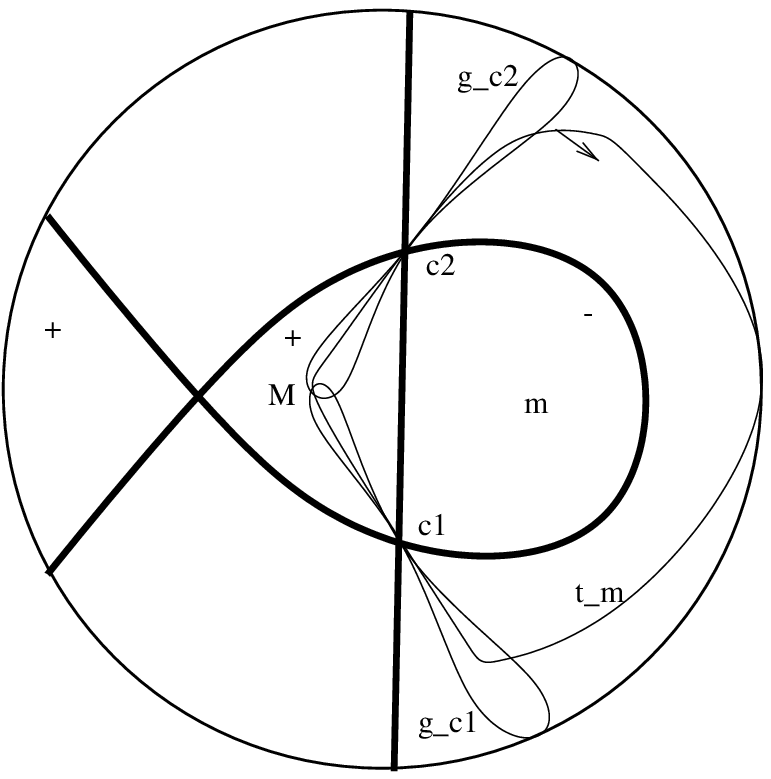}}
\centerline{Fig. $8.$ The vanishing cycle $\delta_M$ for a minimum.}
\endinsert

The link of a divide $P$ is naturally oriented by the following recipe. Let
$\gamma:]0,1[\to D^2$ be a local regular 
parametrization of $P$. The orientation of $L(P)$ is such that the map
$t\in ]0,1[ \mapsto (\gamma(t), \lambda(t)\dot{\gamma}(t))\in L(P)$ is
oriented. Here $\lambda(t)$ is a positive scalar function, which ensures that
the map takes its values in $L(P)$. For a connected divide $P$ we orient
its fiber surface $F_P$ such that the oriented boundary of 
$(F_P \cup L(P),L(P))$ coincides with the orientation of $L(P)$. Vanishing
cycles $\delta_c$ where $c$ is a critical point of $f_P$ do not carry
a natural orientation, since the third power of the geometric monodromy
of the singularity $\{x^3-y^2=0\}$ reverses the orientations. The main use
of a vanishing cycle $\delta_c$ in this paper is through the associated
right Dehn twist $\Delta_c$ of $F_P$. This use 
does not require an orientation for
the cycles $\delta_c$, but requires an orientation of the surface $F_P$. 
Moreover, we orient the tangent space 
$TD=D\times \bR^2$ and its unit sphere $S^3$ as boundary of its unit ball,
such that linking the numbers of links the $L(P_1)$ and $L(P_2)$ are positive 
for generic
pairs of
divides $P_1$ and $P_2$. In fact, the
orientation $TD^2$ is opposite to the orientation as tangent space. 
We have $Lk_{S^3}(L(P_1),L(P_2))=\#(P_1\cap P_2)$, a fact which was 
already used in the previous section.

\bigskip
\S 4. {\bf Reduction cycles and reduction tori}\par

We consider a divide of the form $Q:=P_{p,q}*P$, where the divide $P$ is
given by an immersion $\gamma:[-1,1] \to D^2$. So,
the divide $Q$ is the image of 
$N\gamma \circ T_{p,q}$. If we change the immersion
by a repara\-metrization $\gamma_1:=\gamma \circ \phi$, where $\phi$ is an
oriented diffeomorphism of $[-1,1]$, the divide 
$Q_1:= N\gamma_1 \circ \phi([-1,1])$ is isotopic to 
$Q$ by an transversal isotopy, which does not change 
the type of its knot. By choosing $\phi$ appropriately and $\eta$ small, 
one can achieve that to each double
point of $P$ correspond $p^2$ double points of the divide $Q_1$, 
which look like the intersection points of 
a system of $p$ almost parallel lines with an other system 
of $p$ almost parallel lines, see Fig. $2$, where $p=2$, Fig. $3$, 
where $p=3$ and Fig. $9$. 

\midinsert
\cline{\epsffile {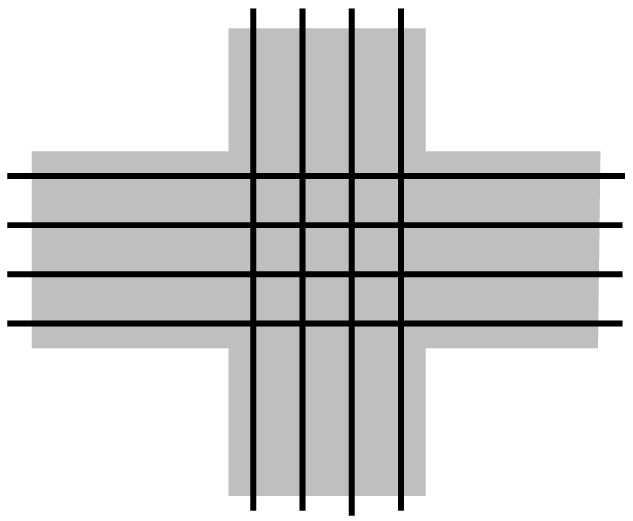}}
\centerline{Fig. $9.$ Manhattan: crossing of the box with $4$ by $4$ strands.}
\endinsert

We may assume, that the divide $Q$ for each double point of $P$
already has a grid of $p^2$ intersections. 

We will construct reduction curves for the monodromy of the knot $L(Q)$
by the method of [AC1]. The reduction curves of the cabling $P_{p,q}*_{\eta}P$
are the intersection of the fiber $F_{Q}$ over 
$1 \in S^1$ of the fibration on
the complement of the knot $L(Q)$ with the boundary of a regular tubular 
neighborhood $U$ of the closed tubular neighborhood $V$ of the knot $L(P)$
forwhich $L(Q) \subset \partial{V}$ holds. The intersection 
$(F_{Q} \cap \partial{U}) \subset F_{Q}$ is indeed 
a system of reduction curves
provided that the torus $\partial{U}$  is transversal to the fibration
of the knot $L(Q)$.

Assume that the tubular neighborhood $V$ was constructed with
the field $\Phi_{p,q}$ and a particular value of the parameter $\eta$.
The same field of sectors, but a slightly bigger parameter value
$\eta'$ yields a tubular neighborhood $U$ of $V$ in $S^3$. The construction
of the fibration will be done as in [AC5]. The main choice for the 
construction of the fibration for the knot $L(Q)$ is a morse function 
$f_Q:D^2 \to \bR$ with $f_Q^{-1}(0)=Q$. For 
our purpose here, where we must achieve the above
transversallity, we will choose $f_Q$ as follows.
First, after applying a  regular transversal small 
isotopy, we may assume that the divide $P$  
has perpendicular rectilinear crossings. Next, we 
consider a morse function $f_P:D^2 \to \bR$ for 
the divide $P$, that is euclidian near its crossings. Let the fibration 
on the complement of the knot $L(P)$ be
$\pi_{P,\eta}:S^3 \setminus L(P) \to S^1$ 
where
$\pi_{P}(x,u):=\theta_{P}(x,u)/ |\theta_{P}(x,u)|$
and
$$
\theta_{P,\lambda}(x,u):=f_P(x)+i \lambda^{-1}\, df_P(x)(u)-
{1\over 2}\lambda^{-2}\chi(x)H_{f_P}(x)(u,u).
$$
The function $\chi:D^2 \to \bR$ is a bump function at the crossing 
points of $P$ and $\lambda$ is a big real parameter.
We now choose  a small positive real number $v$, such that 
$\{x \in D^2 \mid |f_P(x)|\leq v\}$ is a regular tubular of $P$, that
meets each component of $\{x\in D^2 \mid \chi(x)=1\}$. 
Next 
we choose $\eta'>0$ such the corners of $N_{\eta'}\gamma(B)$ are in
$\{|f_P(x)|=v\}$ i.e. $\eta'^2=v$.  We constuct the torus knot $L(Q)$ with
$Q:=P_{p,q}*_{\eta}P$ where $0 < \eta < \eta'$. Since 
$Q \subset \{|f_P(x)| < v\}$ holds, we can construct 
a morse function
$f_Q:D^2 \to \bR$ for the divide $Q$, such 
that on $\{|f_P(x)| \geq v\}$ the function $f_Q$
is constant along the level sets of $f_P$. 

The following 
theorem follows directly from 
Lemme $2$, page $153$ in [AC1] and the above construction.

\proclaim{Theorem $1$} The torus $\partial{\Phi_{\eta',p,q}}\gamma$ 
is transversal to the fibration
deduced from $f_Q$ on the complement of the knot $L(Q)$. The intersection 
$$
\partial{\Phi_{\eta',p,q}}\gamma \cap F_{Q}
$$ 
is a system of $p$ closed  curves on the fiber $F_{Q}$, which 
is a reduction of the monodromy of $L(Q)$.
\endproclaim

With a few examples, we now explain how to depict in the Milnor fiber
a distinguished system of vanishing cycles and the
reduction curves for the monodromy of a singularity,
for which a divide of the form $Q=P_{p,q}*P$ is given. 

The fiber $F_{Q}$ with a distinguished system of
vanishing cycles is already constructed in Section $3$. 

For a $(p,q)$ cabling the reduction 
system consists of $p$ simply closed curves on the fiber $F_{Q}$. Each 
of them cuts out from $F_{Q}$ a surface diffeomorphic to the fiber 
$F_{P}$ of the divide $P$. The $p$ copies of $F_{P}$ in $F_{Q}$ 
are cyclicly permuted by the monodromy $T_Q$.

One of those copies can be visualized more easily as follows.
Let $\{x \in D^2 \mid f_P(x) \geq v \}$. For each double point $c$ of $P$
we connect the two components of $\{x \in D^2 \mid f_P(x) \geq v \}$,
that are incident with the double point $c$, by a special bridge  
which projects diagonally through
the Manhattan part  of the divide $Q$,
that corresponds to $c$. The
projection of the bridge is  a twisted strip $S_c$ in $D^2$, that 
realizes a boundary connected sum of
the $P_+$-components. The twist points of the strip $S_c$ are precisely
the critical points of $f_Q$, that lie on the diagonal. The boundary
of $S_c$ consists of two smooth curves, that intersect
each other transversally and that also intersect the divide $Q$ 
transversally. 

\midinsert
\cline{\epsffile {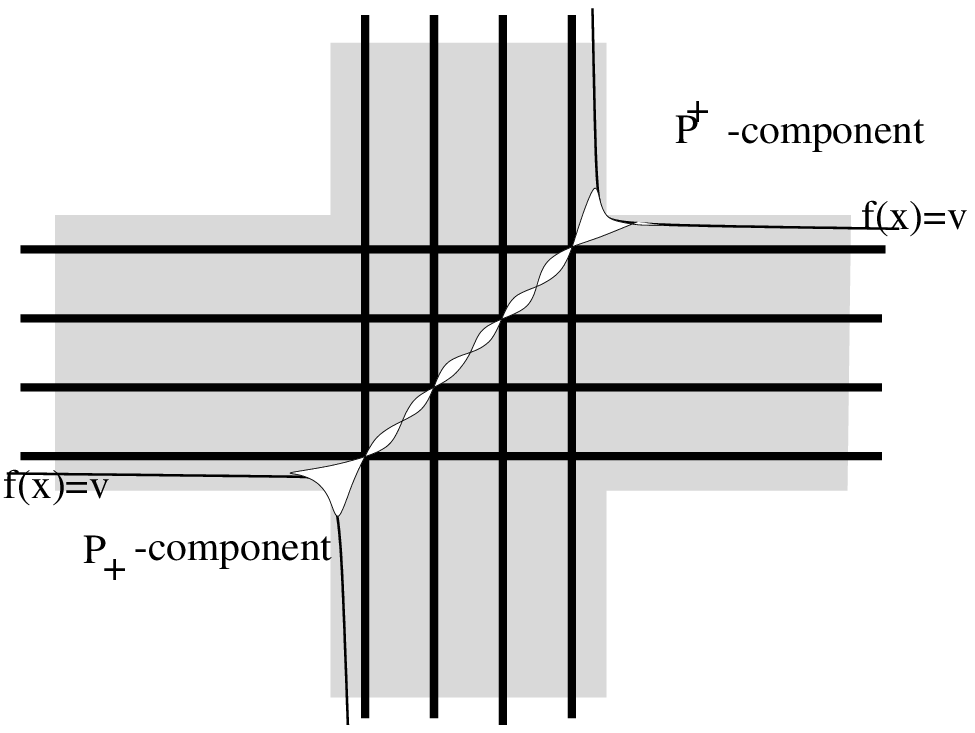}}
\centerline{Fig. $10.$ Bridge through Manhattan.}
\endinsert

Let $C$ be the union
of the projections $S_c$ of the 
bridges with $\{x \in D^2 \mid f_P(x) \geq v \}$, see Fig. $10$. 
In Fig. $11$ we have zoomed
out one block to  show more details.
The copy $F_{P,Q}$ 
of the fiber of the knot $L(P)$ is the closure  in the fiber $F_Q$
of the knot $L(Q)$ of the set
$$
\{(x,u) \in TD^2 \mid x \in C, (df_Q)_x\not=0, (df_Q)_x(u)=0\}\cap S^3.
$$

The first
reduction curve   is the 
boundary of the surface $R:=\partial{F_{P,Q}}$ 
of the surface
$F_{P,Q}\subset F_Q$. 

The reduction system
is the orbit $\{R, T_Q(R), T_Q^2(R), \dots \}$ under 
the monodromy $T_Q$ of the singularity with divide $Q$.

\midinsert
\cline{\epsffile {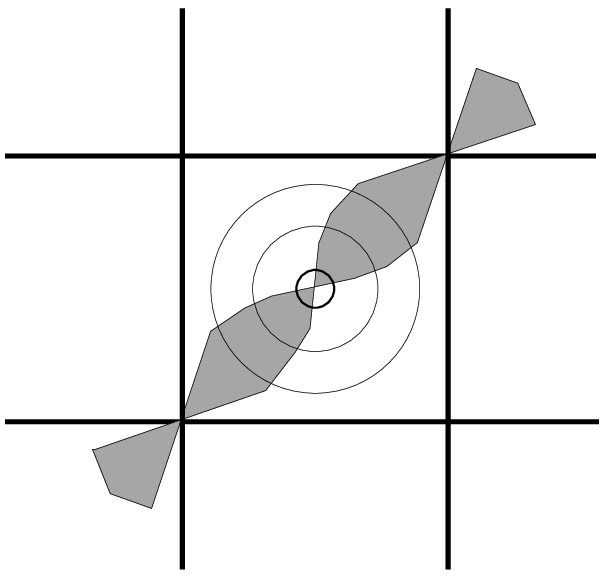}}
\centerline{Fig. $11.$ Bridge through a block of Manhattan.}
\endinsert

Our first example is the singularity with two essential Puiseux pairs
$(x^3-y^2)^2-4x^5y-x^7$, whose link is a two stage iterated torus knot. Its
divide $Q=P_{2,9}*P_{2,3}$, (see Fig. $2$), has two $P_+$-components, where
$P=P_{2,3}$, 
(see Fig $12$, where the projection of the reduction curve $R$ is
drawn).  In this case Manhattan consists of one block. 
The reduction curve $R$ is the pre-image in the fiber 
$F_Q$ of its projection $\proj(R) \subset D^2$ under
the map $(x,u) \mapsto x$ is drawn in Fig. $12$. That means $R$ 
is the closure in $F_Q$ of the set
$$
\{(x,u) \in TD^2 \mid x \in \proj(R), (df_Q)_x\not=0,(df_Q)_x(u)=0, ||x||^2+||u||^2=1\}.
$$
The reduction system is $\{R, T_Q(R)\}$.

\midinsert
\cline{\epsffile {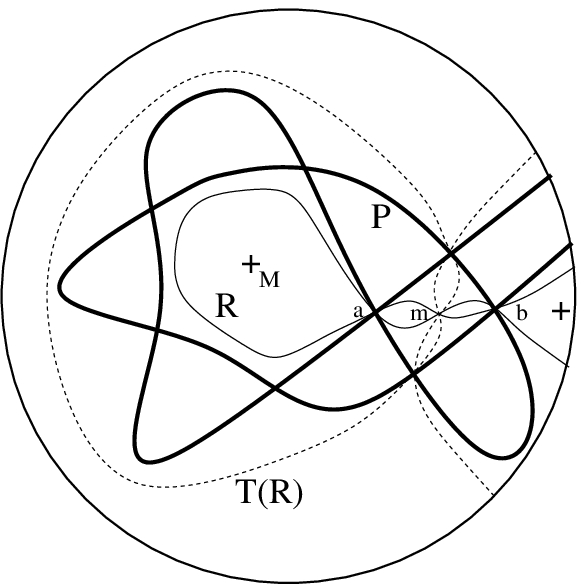}}
\centerline{Fig. $12.$ Divide $Q=P_{2,9}*P_{2,3}$ with reduction curves $R$ and $T(R)$ (dotted).}
\endinsert

The curve $R$ is homologically trivial in $F_Q$. It turns out
that the power $T_Q^{156}$ of the monodromy is the composition of the
right Dehn twists, whose core curves are  
$\{R, T_Q(R)\}$. The power
$T_Q^{156}$ is a product of $2496=16 \times 156$ Dehn 
twists, since $T_Q$ is the
product of those Dehn twists whose core 
curves are the system of distinguished
quadratic vanishing cycles of the real morsification with divide $Q$. 
It   turns out that the expression as product
of Dehn twists is far from being as short as possible. In fact, the right 
Dehn
twist $\Delta_R$ with core curve $R$ can be 
written as a product of $36$ right Dehn twists
that have core curves coming from the morsification with divide $Q$. 
More precisely, the Dehn twist $\Delta_R$ factors as
$$
\Delta_R=(\Delta_M \circ \Delta_b \circ \Delta_m \circ \Delta_a \circ \Delta_m^{-1} \circ \Delta_b^{-1})^6,
$$
The factors are right Dehn twists whose core curves are 
among the quadratic 
vanishing cycles 
$\delta_m,\delta_a,\delta_M,\delta_b$ of the 
divide $Q$ as indicated in Fig. $11$, 
$\delta_M$ is the
vanishing cycle of a $P_+$-region, $\delta_m$ of the 
maximum of Manhattan, and $\delta_a,\delta_b$
of street corners of Manhattan. It follows 
that $T_Q^{156}$ can also be written 
as a composition of $72$ Dehn twists with core curves 
among the vanishing cycles of  the divide $Q$. 
The 
composition $\Delta_b \circ \Delta_m \circ \Delta_a \circ \Delta_m^{-1} \circ \Delta_b^{-1}$ is the Dehn 
twist with core curve $\bar{a}:=\Delta_b(\Delta_m(a))$.

\midinsert
\cline{\epsffile {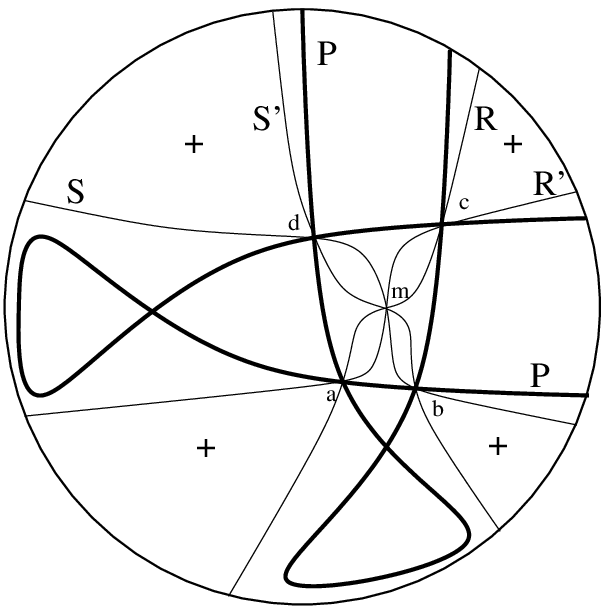}}
\centerline{Fig. $13.$ Divide $P$ for $(x^3-y^2)(y^3-x^2)$ with reduction system  $R \cup S$.}
\endinsert

The reduction curve $R$ cuts off from $F_Q$ a piece $F_{P,Q}$ of genus 
one (see also
Fig. $1$ on page $159$ of [AC1]), and the Dehn
twists $\Delta_M$ and $\Delta_{\bar{a}}$ act only on this piece, since
the curves $\delta_M$ and $\bar{a}$ lie entirely in this piece; in
this piece, that is a copy of the fiber $F_P$, they generate
the geometric monodromy group of the accompanying
singularity $x^3-y^2=0$ with divide $P_{2,3}$.

Our second example is the singularity with two branches $(x^3-y^2)(y^3-x^2)$.
Its homological monodromy is of infinite order [AC1].
Each branch is a torus knot. Again Manhattan consists of 
one block. In Fig. $13$
we have drawn the projections of the curves 
$R,R'$ and $S,S'$, that together are the boundary components of the 
two diagonals through Manhattan. In this case the curves $R$ and $R'$ are
isotopic to each other, as are the curves $S$ and $S'$. A complete reduction
system for the geometric monodromy is the system $\{R,S\}$.
Each component
of this system carries a non-trivial homology class. The isotopy classes
of the curves $R$ and $S$ are permuted by the monodromy $T_P$, hence the system
$\{R,S\}$ is invariant under the monodromy. 

Let $h$ be the 
action of $T_P$  on
the homology $H_1(F_P,\bZ)$ of the the fiber $F_P$. Let
$\delta_a, \delta_b, \delta_c, \delta_d$ be the 
vanishing cycles of the double points, that are the
corners of Manhattan of $P$, let $\delta_m$ be the vanishing cycle of the
maximum in the center of Manhattan.

If one chooses the orientations appropriately, on has 
$$
[R]=[\delta_a]+[\delta_m]+[\delta_c],\, [S]=[\delta_b]+[\delta_m]+[\delta_d],\, h([R])=-[S],\, h([S])=-[R],$$ 
hence also $h([R]-[S])=[R]-[S]$. Let $[k]$ be any cycle on $F_P$, that
is carried by a simple oriented curve $k$ and intersects the curves $R$ and $S$
each transversally in one point. One has $h^{10}([k])=[k]\pm([R]+[S])$,
which shows that the homological monodromy $h$ is not of finite order. We
have drawn in Fig. $13$ the oriented projection of such a cycle $k$, that
intersects the curves $A$ and $B$. 
The curve $A$ is halfway in between the curves $R$ and
$R'$ on the cylinder they cut out. Let $B$ be the curve halfway in between
$S$ and $S'$. The curves $A$ and $B$ are the reduction curves of Fig. $14$
on page $167$ of [AC1]. The reduction system $A,B$ is much easier to
draw, see Fig. $14$, where are drawn the projections in $D^2$. The projections
meet transversally at the maximum in Manhattan of $f_P$.
The curve
$\delta_m$ intersects 
transversally in two points each curve $R$ and $S$. One has 
$h^{10}([\delta_m])=[\delta_m]\pm 2([R]+[S])$.

\midinsert
\cline{\epsffile {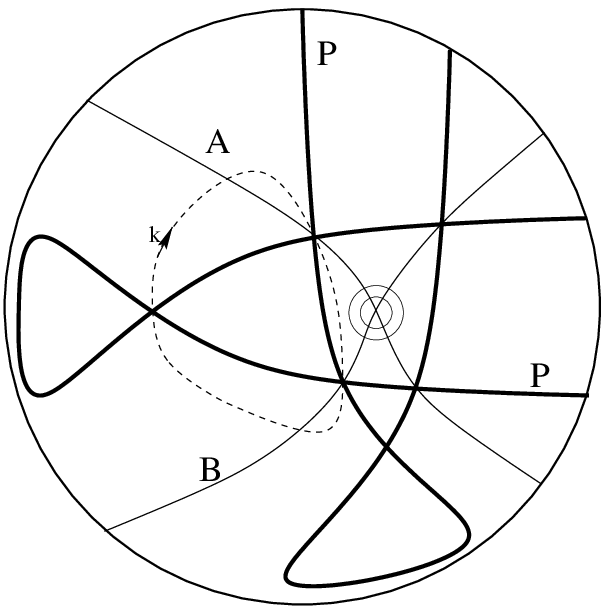}}
\centerline{Fig. $14.$ Divide $P$ for $(x^3-y^2)(y^3-x^2)$
with reduction system  $A \cup B$.}
\endinsert

The power $T_P^{10}$
of the geometric 
monodromy, that is a word of length 
$110$ in the Dehn twists of the divide $P$,
is equal to the composition of those right Dehn
twists, whose core curves are  $R$ and $S$. So, the power $T_P^{10}$
also can be written as  the much shorter word
$$
\Delta_c \circ \Delta_m \circ \Delta_a \circ \Delta_m^{-1}\circ 
\Delta_c^{-1}\circ 
\Delta_d \circ \Delta_m \circ \Delta_b \circ \Delta_m^{-1} \circ \Delta_d^{-1}.
$$

\bigskip
{\bf Remark.} The curves $A$ and
$\Delta_m(\Delta_c(\delta_a))$ are isotopic, where 
$\Delta_m$ and $\Delta_c$ are the right Dehn
twists with core curves $\delta_m$ and $\delta_c$. It 
follows that the reduction system
$A,B$ consists  of quadratic vanishing cycles
of the singularity of $\{(x^3-y^2)(y^3-x^2)=0\}$ with two branches. 
In contrast, a  reduction curve of a 
singularity with only one branch can not be   a quadratic vanishing cycle,
since all reduction curves are zero in the homology.  
\bigskip

\S 5. {\bf Geometric monodromy group and reduction system.}

Let the polynomial $f_{(a,b)}$ be an equation for an irreducible
plane curve singularity with $n$ 
essential Puiseux pairs $(a_i,b_i)_{1\leq i \leq n}$. The number of 
simply closed curves contained in a complete reduction
system $R$ for the monodromy of $f$ is 
$$
a_na_{n-1}\cdot\, \dots \,\cdot a_2+a_{n-1}a_{n-2}\cdot \, \dots \,\cdot a_2+ \dots +a_3a_2+a_2.
$$
Let $\Gamma_{f,red}$ be the subgroup of the 
geometric monodromy group
of $\Gamma_f$ of $f$ of those elements 
$\gamma \in \Gamma_f$ that up to isotopy fix
each component of $R$. Let $\Gamma_{f,red}^0$ be the subgroup of $\Gamma$
which is generated by the Dehn twist whose core curves are
quadratic vanishing cycles and do not intersect any component of $R$.
Obviously, one has $\Gamma_{f,red}^0 \subset \Gamma_{f,red}$, but 
we do not know if this inclusion is strict. A component of $F \setminus R$ is
called a top-component if its closure in $F$ meets only one
component of $R$. Let $\Gamma_{f,top}$ be 
the subgroup of $\Gamma_{f}$ of those monodromy transformations, 
which induce the identity in each component of $F \setminus R$ that 
is not a top-component. 
Let $\Gamma_{f,top}^0$ be the intersection 
$\Gamma_{f,top}\cap \Gamma_{f,red}^0$. We have

\proclaim{Theorem $2$} Let $f=f_{(a,b)}$ be an irreducible singularity with
$n \geq 2$ essential Puiseux pairs $(a_i,b_i)_{1\leq i \leq n}$. 
Let $g=f_{(a',b')}$ be a singularity with 
the $n-1$ essential Puiseux pairs
$(a',b')=(a_i,b_i)_{1\leq i \leq n-1}$. The group 
$\Gamma_f$ contains the product
of $a_n$ copies of the group $\Gamma_g$.
\endproclaim

\proclaim{Theorem $3$} Let $f_{(a,b)}$ be an irreducible singularity with
$n \geq 2$ essential Pui\-seux pairs. The group $\Gamma_{f,top}^0$ is isomorphic
to the product of $a_na_{n-1}\cdot\, \dots \,\cdot a_2$ copies of the
geometric monodromy group of the
singularity $y^{a_1}-x^{b_1}=0$.
\endproclaim

{\bf Proof of Th. $2$.} Let $P$ be the divide 
$P_{a_{n-1},b'_{n-1}}* \dots *P_{a_2,b'_2}*P_{a_1,b_1}$ for the singularity
of $g$ and let $Q=P_{a_n,b'_n}*P$ be the divide for the singularity of $f$.
A copy $F_{P,Q}$ of the fiber $F_P$ is constructed as a subset of the fiber
$F_Q$. Remember, that $F_P$ is obtained by connecting with strips
the sets $\{(x,u)\in TD^2 \mid f_P(x)>0, (df_P)_x(u)=0\}$, where $f_P:D^2 \to \bR$
is a morse function for the divide $P$. For each double point of $P$
there are two connecting strips. To each $+$-component of $P$ corresponds a
$P_+$-component of $Q$ with the same topology  
and to each double point of $P$ corresponds a Manhattan
grid of $Q$, in which we have drawn diagonally the projection of the strips
that connect $\{(x,u)\in TD^2 \mid x \in Q_{P,+}, (df_Q)_x(u)=0\}$. Here,
$Q_{P,+}$ denotes the union of the $P_+$-components  of the 
complement of the divide $Q$. From the divide $P$ is deduced a distinguished
base of quadratic 
vanishing cycles for the singularity of $f$. Let 
$B_P$ be the union of the curves of this base.  This base can be drawn
on the fiber $F_P$, see Section $3$. 

In order to prove the theorem, we will 
construct inside $F_{P,Q}$ a system
of simply closed curves with union $B_{P,Q}$, each of them being a 
quadratic vanishing cycle for
the singularity $g$, such that the pairs $(F_P,B_P)$ and 
$(F_{P,Q},B_{P,Q})$ are diffeomorphic. This finishes the proof,
since the Dehn twist, whose cores are the quadratic vanishing cycles
of $B_{P,Q}$, generate a copy of $\Gamma_g$ in $\Gamma_f$. By acting with the 
geometric monodromy $T$ of the singularity $f$, one obtains $a_n$ commuting
copies of $\Gamma_g$ in $\Gamma_f$.

To each $+$-region of $P$ corresponds one $P_+$-region of $Q$. The maximum
of $f_P$, say at $M$ in the region, is also a maximum of $f_Q$. The quadratic
vanishing cycle $\delta_M:=\{(M,u)\in TD^2 \mid ||M||^2+||u||^2=1\}$ of $F_Q$ 
lies in $F_P$
and also in $F_{P,Q}$.  For each double point $c$ of $P$ the quadratic
vanishing cycle $\delta_c \subset F_P$ projects in $D^2$ to a tear splicing
the
gradient line $g_c$ of $f_P$ through $c$. The endpoints of $g_c$ 
are maxima of $f_P$ or points on $\partial{D^2}$. The function $f_Q$ 
has exactly one gradient line $g_{Q,c}$ that has the same endpoints
as $g_c$ and coincides with $g_c$ 
in a neighborhood of the common endpoints. The gradient line $g_{Q,c}$
runs along a diagonal through the Manhattan grid corresponding to $c$. Let
$g_{Q,c}$ be the simply closed curve on $F_Q$, that projects to a tear
$t_{Q,c}$ equal to $g_{Q,c}$, except  above a neighborhood of its endpoints
where $t_{Q,c}$ equals $t_c$. 
We remark that $\delta_{Q,c}$ is a cycle in $F_{P,Q}$. Let
$c_1, \dots ,c_p$ be the $p:=a_n$ double points of $Q$ that occur along the 
$g_{Q,c}$ and let $M_2, \dots , M_p$ along $g_{Q,c}$ be 
the maxima.
Let $\delta_{Q,c_1}$ be the quadratic vanishing cycle of the singularity $f$
that corresponds to $c_1$. One verifies that the cycles
$\delta_{Q,c}$ and 
$$
\Delta_{c_p}\circ \Delta_{M_p} \circ \dots \circ \Delta_{c_2}\circ \Delta_{M_2}(\delta_{Q,c_1})
$$
are isotopic. Here $\Delta_{c_i}$ or $\Delta_{M_i}$ stands for 
the right Dehn twist of $F_Q$ whose
core curve is the quadratic vanishing cycle 
$\delta_{c_i}$ or $\delta_{M_i}$
of the singularity $f$. 
Hence $\delta_{Q,c} \subset F_{P,Q}$ is a quadratic
vanishing cycle for the singularity $f$. So far, we have constructed for
each maximum and for each saddle point of $f_P$ a simply closed curve on
$F_{P,Q}$ that is a quadratic vanishing cycle of the singularity $f$. These
cycles intersect on $F_{P,Q}$ as do the corresponding quadratic vanishing 
cycles of the singularity $g$ on $F_P$. 

\midinsert
\cline{\epsffile {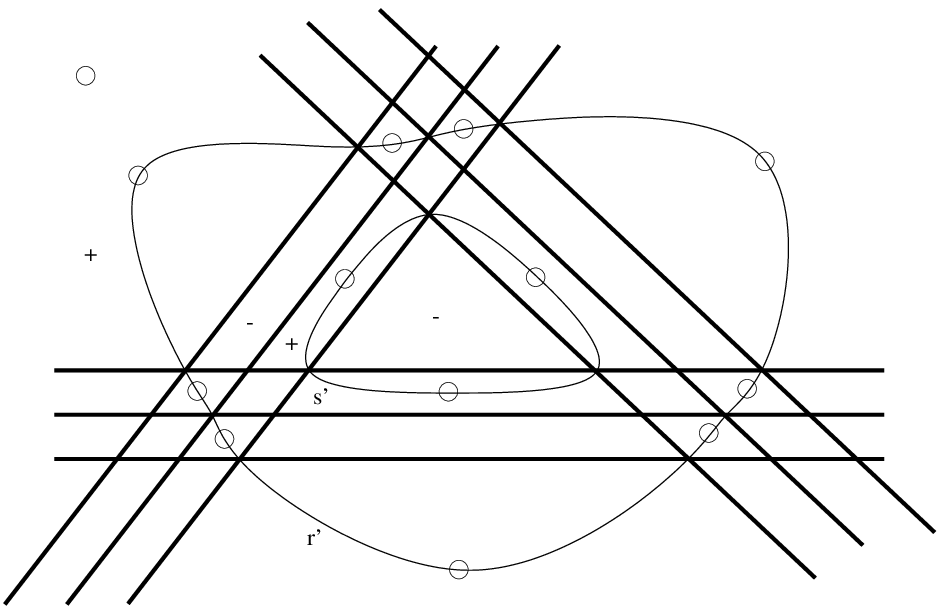}}
\centerline{Fig. $15.$ Vanishing cycle $s$ on $F_Q$ from a minimum of 
$f_P$ and cycle $r$ on $F_{P,Q}$.}
\endinsert

We now wish to construct for each minimum of $f_P$ a vanishing 
cycle on $F_{P,Q}$. We have to handle two cases: $p$ odd, see Fig. $15$, and
$p$ even, see Fig. $16$.

If $p$ is odd, a minimum $m$ of $f_P$ will 
also be a minimum of $f_Q$. Let  $\delta_{Q,m}$ 
be the vanishing cycle on $F_Q$ corresponding to $m$, see Fig. $15$. The 
projection of $\delta_{Q,m}$ into $D^2$ is a smooth simply closed 
curve $s'$ transversal to $Q$, that 
surrounds the $-$ region of $m$ through
the its neighboring $+$ regions of $Q$. One needs to take care 
that in each neighboring
$+$ component the projection runs through the maximum of $f_Q$ in that region.
The points of $s$ correspond to pairs $(x,u)$ with $x\in s'$ and $u$ pointing
inwards to $m$. Let $r$ be a simply closed cycle on $F_{P,Q}$ 
that projects into
$D^2$ upon the curve $r'$, which now surrounds the $-$ region of $m$ 
through
the $P_+$-components of $Q$, see Fig. $15$. In the 
Manhattan grids $r'$ ist just
a diagonal, again $r'$ runs through the maxima of the regions or 
touches $\partial{D}$. On $r$ we 
only allow pairs $(x,u)$ where $u$ points inwards to $m$. It is clear that 
the cycle $r$ on $F_{P,Q}$ intersects the cycles of the previous
construction as the vanishing cycle to the minimum of $f_P$
intersects the vanishing cycles of the critical points of $f_P$.
It remains however to check that the cycle $r$ is a quadratic
vanishing cycle of the singularity of $g$. By 
applying to $r \subset F_Q$ 
the Dehn twist corresponding to the
critical points of $f_Q$ that are in between the curves $r'$ and $s'$,
one can transform the isotopy class of the curve $\delta_{Q,m}$ 
to the class of the curve
$r$. This proves that $r$ is indeed a quadratic vanishing cycle of the
singularity of $g$.

\midinsert
\cline{\epsffile {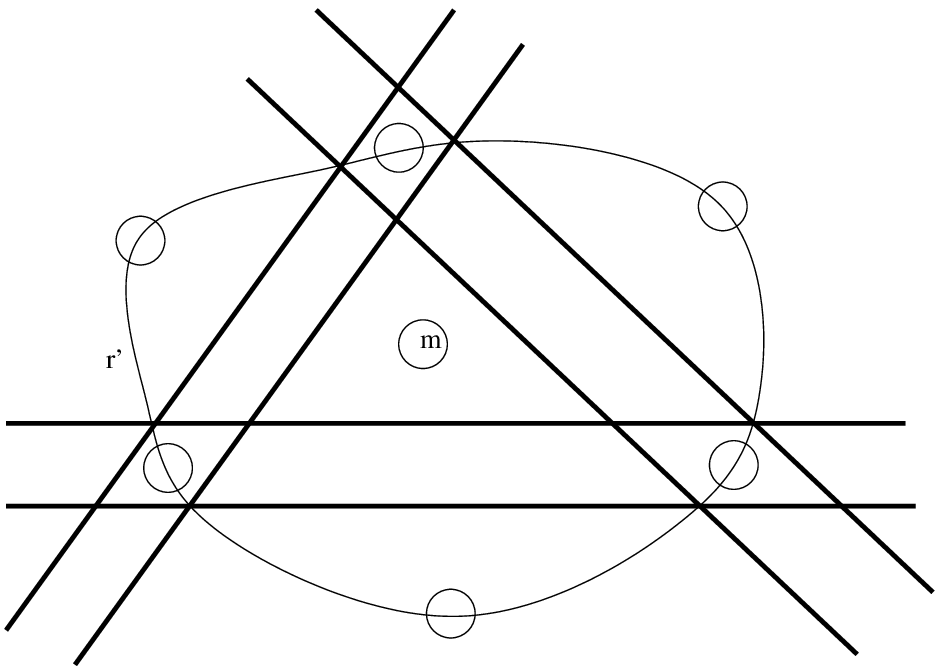}}
\centerline{Fig. $16.$ Vanishing cycle $\delta_{Q,m}$ on $F_Q$ from a minimum of
$f_P$ and cycle $r$ on $F_{P,Q}$.}
\endinsert 

If $p$ is even, a minimum $m$ of $f_P$ will
be a maximum of $f_Q$. Let  $\delta_{Q,M}$
be the vanishing cycle on $F_Q$ corresponding to maximum $M:=m$, 
see Fig. $16$. Its 
projection into $D^2$ is the point $M:=m$. Let $r$ be a simply 
closed cycle on $F_{P,Q}$
that projects into
$D^2$ upon the curve $r'$ which now surrounds the $-$ region of $M:=m$
through
the $P_+$-regions of $Q$, see Fig. $16$. In the Manhattan grids $r'$ ist just
a diagonal, again $r'$ runs through the maxima of the regions. On $r$ we
only allow pairs $(x,u)$ where $u$ points inwards to $m$. It is clear that
the cycle $r$ on $F_{P,Q}$ intersects the cycles of the previous
construction as the vanishing cycle to the minimum of $f_P$
intersects the vanishing cycles of the critical points of $f_P$.
By
applying to $\delta_{Q,M} \subset F_Q$
the Dehn twist corresponding to the
critical points of $f_Q$ that are in between the curve $r'$ and 
the point $M:=m$,
one can transform the isotopy class of the curve $\delta_{Q,M}$ to the class of the curve
$r$, and proves that $r$ is indeed a quadratic vanishing cycle of the
singularity of $g$. As explained, this terminates the proof.

\bigskip

{\bf Proof of Th. $3$.} The proof of Theorem $2$ constructs a copy 
$\Gamma_{P,Q}$ of the
mono\-dromy group $\Gamma_f$ of the singulatity $f$ as subgroup 
in the monodromy group $\Gamma_g$ of the
singularity $g$. This copy acts with support in a copy $F_{P,Q}$ 
of the fiber $F_P$. The the first $a_n-1$ iterates of the 
monodromy $T_Q$ of the singularity $g$  constructs $a_n$ copies of $F_P$
in $F_Q$. By conjugation with $T_Q$ one gets $a_n$ copies from $\Gamma_{P,Q}$.
We end the proof by repeating this argument. One 
gets $a_na_{n-1} \cdots a_2$ commuting copies 
of the geometric monodromy group of the 
singularity $x^{b_1}-y^{a_1}$ in $\Gamma_g$.

\bigskip

{\bf Problem.} We like to state the problem of presenting
the geometric
monodromy group of plane curve  singularities with generators
and relations. It 
would be particulary nice to express
the presentation in terms of a divide of the singularity. The same
problem can also be stated for the homological monodromy group of
plane curve singularities, but I think that the problem for the
geome\-tric monodromy group is more tractable, since all reduction curves
can be 
taken into account. The theorems $2$ and $3$ are possibly first steps towards
a solution of this problem. However, an important missing 
piece in this program is a presentation
with generators and relations of the geometric monodromy group of the
singularities $y^p-x^q=0$ for $3\leq p \leq q,\, 7 \leq p+q$. The fundamental
group of the complement of the discriminant in the unfolding of
the singularity $y^2-x^q=0$ is the braid group $B_{q-1}$. Bernard Perron and 
Jean-Pierre Vannier have proved for the singularities $y^2-x^q=0$ that
the geometric monodromy group is a faithful image of the braid group $B_{q-1}$
and that a similar result holds for the singularities $x(y^2-x^q)=0$ [P-V].
The fundamental
group of the complement of the discriminant in the unfolding of
the singularity $y^3-x^6=0$ is the Artin $AE_6$ group 
of the Dynkin diagramm $E_6$.
Bronek Wajnryb has proved that the geometric 
monodromy representation of $AE_6$ into the mapping class group 
of the Milnor fiber of the singularity $y^3-x^6=0$ is not faithful.

\goodbreak
\par

\bigskip
\noindent 
Universit\"at Basel \par
\noindent 
Rheinsprung 21 \par
\noindent 
CH-4051  Basel.

\bye